\documentstyle[11pt]{amsart}

\newtheorem{prop}{Proposition}
\newtheorem{theo}{Theorem}
\newtheorem{Lemma}{Lemma}
\newtheorem{cor}{Corollary}

\def\Bbb{\bf}

\newcommand{\na}{\nabla}

\newcommand{\Om}{\Omega}
\newcommand{\si}{\sigma}

\newcommand{\La}{\Lambda}

\newcommand{\miu}{\mu}

\newcommand{\ka}{K{\"a}hler }

\title[Almost K{\"a}hler 4-manifolds with $J$-invariant Ricci tensor ....]
{Almost K{\"a}hler 4-manifolds with $J$-invariant Ricci tensor and special
Weyl tensor}
\author{Vestislav Apostolov and Tedi Dr\u{a}ghici}
\address{Vestislav Apostolov \\ Institute of Mathematics and Informatics\\
Acad. G.Bonchev St.\\ Bl. 8\\ 1113 Sofia, Bulgaria} 
\address[Current address]{Centre de
Math{\'e}matiques\\ Ecole Polytechnique \\ U.M.R. 7640 du C.N.R.S. \\
F-91128 Palaiseau \\ France}
\email{apostolo@@math.polytechnique.fr}
\address{Tedi Dr\u{a}ghici\\ Department of Mathematics\\ Northeastern
Illinois University\\
Chicago IL 60625-4699\\ USA}
\email{TC-Draghici@@neiu.edu}

\thanks{The first-named author is partially supported by a grant of the 
EPDI/IHES}

\begin{document}
\date{}
\maketitle

\section{Introduction}

An {\it almost K{\"a}hler structure} on a manifold $M^{2n}$ is a triple
$(g, J, \Om)$
of a Riemannian metric $g$, an orthogonal almost complex structure $J$, and a
symplectic form $\Om$
satisfying the compatibility relation
\begin{equation}
\Om(X,Y) = g(JX, Y),
\end{equation}
\noindent
for any tangent vectors $X,Y$ to $M$. If the almost complex structure $J$
is
integrable we obtain a {\it K{\"a}hler structure}. Many efforts have been
done in
the direction of finding curvature conditions on the metric which insure
the integrability of the almost complex structure. A famous conjecture of
Goldberg \cite{Go} states that a compact almost K{\"a}hler, Einstein
manifold
is in fact K{\"a}hler. Important progress was made by K. Sekigawa who
proved
that the conjecture
is true if the scalar curvature is non-negative \cite{Se2}.
The case of negative
scalar curvature is still wide open, although recently there has been some
significant progress in dimension 4, see \cite{Ar2,Ar3,Ar,MOS,NuP,OS2}.
It is now
known that if the conjecture
turns out to be true, the compactness should play an essential role.
Nurowski and Przanowski \cite{NuP} constructed a 4-dimensional, local
example of Einstein, strictly almost \ka  manifold; this was
generalized by K.Tod (see \cite{Ar3,Ar}) to give a family of such
examples. It is interesting to 
remark that the structure of the Weyl tensor of all these examples is
unexpectedly special. An important recent result of J. Armstrong
\cite{Ar3,Ar} 
states that any 4-dimensional almost K{\"a}hler, non-K{\"a}hler, Einstein 
manifold is obtained by Tod's construction, provided that the
K{\"a}hler form is an eigenform of the Weyl tensor.
In the compact case, on the other hand,
some of the positive partial results
on the conjecture in dimension 4 have been obtained exactly by imposing
some additional assumptions on the structure of the Weyl tensor 
(\cite{Ar2,Ar3,Ar,OS2}).

The recently discovered Seiberg-Witten invariants could have an
impact towards a complete answer to the Goldberg conjecture in dimension
4.
These are invariants of smooth, oriented, compact 4-dimensional manifolds.
The works of Taubes
\cite{Ta2,Ta1,Ta3} and others showed the invariants to be particularly
interesting for
symplectic 4-manifolds. (See also \cite{Ko3} for a quick introduction to
Seiberg-Witten invariants and some of their applications to symplectic
geometry.)
LeBrun \cite{LeB1}, \cite{LeB2}, \cite{LeB3} and Kotschick
\cite{Kotschick} have also proved several
interesting applications of the Seiberg-Witten theory to the existence of
Einstein metrics in dimension 4, which shed some light on the Goldberg
conjecture as well.

\vspace{0.2cm}
In this paper we will not directly work on the Goldberg conjecture, but
rather
on problems which are parallel to it. The curvature of a K{\"a}hler
metric has strong symmetry properties with 
respect 
to $J$. On the other hand, an arbitrary almost Hermitian (or almost
K{\"a}hler)
metric may not have any of these symmetry properties. A. Gray \cite{Gr}
introduced classes of almost Hermitian manifolds whose curvature
tensor has a certain degree of resemblance to that of a \ka manifold
(see Section 3). The condition that the Ricci tensor is $J$-invariant
could be considered to be the minimal degree of resemblance. 
This is weaker than the Einstein condition and may be even more
natural
in the context of almost Hermitian geometry, as some interesting variational problems on compact symplectic manifolds 
have lead to almost \ka metrics with
$J$-invariant Ricci tensor. It was shown in
\cite{BI}
that they are the critical
points of the Hilbert functional, the integral of the
scalar curvature, restricted to the set of all
compatible metrics to a given symplectic form. Compatible \ka metrics 
provide absolute maxima for the functional in this
setting \cite{Bl}. Blair and Ianus asked if the $J$-invariance of the 
Ricci tensor is a sufficient condition for the integrability of 
an almost \ka structure on a compact manifold \cite{BI}. If true,
this would be a stronger statement than the Goldberg
conjecture. However,  6-dimensional counter-examples to the question
of Blair and Ianus were found by Davidov and Mu\u{s}karov \cite{DM}
(see also \cite{AGI}). Multiplying these by compact K{\"a}hler
manifolds, one 
obtains counter-examples in any dimension $2n, n\ge 3$ .
No example of a strictly almost \ka
structure with $J$-invariant Ricci
tensor is known yet on a compact 4-dimensional manifold.

The examples of \cite{DM} have, in fact, a higher degree
of resemblance to \ka structures than just the $J$-invariant Ricci tensor.
Indeed, they are strictly almost \ka manifolds satisfying the second
curvature condition of Gray (see Section 3 for definition and \cite{Gr}).
In dimension 4, the second curvature
condition of Gray on an almost Hermitian metric just means that the Ricci
tensor and the positive Weyl tensor have the same symmetries as for a K{\"a}hler
metric. 
The integrability of the almost complex structure follows even locally
for an Einstein, almost \ka 4-dimensional manifold which 
satisfies the second curvature condition of Gray
\cite{Ar2,Ar,OS2}. The main result  of this paper is that compact,
4-dimensional, {\it strictly} almost \ka manifolds satisfying the
second curvature condition of Gray could exist only under strong
topological assumptions:  
\begin{theo}\label{th1} 
A compact, 4-dimensional, almost \ka 
manifold  which satisfies the second curvature condition
of Gray  is K{\"a}hler, provided that either the Euler characteristic or the
signature of the underlying manifold does not vanish. 

Moreover, if there exists a compact,
4-dimensional, {\rm non-K{\"a}hler}, almost \ka  
manifold  $(M,g,J)$ satisfying the second curvature condition
of Gray, then $(M,g)$ admits an orthogonal \ka structure ${\bar J}$,
compatible with the reversed orientation of $M$, such that $(M,{\bar
J})$ is a minimal class $VI$ complex surface. 
\end{theo}

\noindent
{\bf Note}: Some time after this paper was completed, the authors
in collaboration with D. Kotschick proved  that the exceptional case of the
Theorem 1 cannot occur  ({\it An integrability theorem for almost
K\"ahler 4-manifolds}, C. R. Acad. Sci. Paris, t. 329, s\'er. I
(1999), 413-418). Thus, any compact almost \ka 4-manifold which satisfies
the second curvature condition of Gray is, in fact, K\"ahler.

\vspace{0.2cm}

\noindent
Note the
contrast with the Hermitian case. For Hermitian surfaces the second
curvature condition of Gray is equivalent to the $J$-invariance of the
Ricci tensor; a number of compact, non-K{\"a}hler
Hermitian surfaces with $J$-invariant Ricci tensor have been
provided on rational surfaces (\cite{AG}).

\vspace{0.2cm}

Another natural problem  which involves the $J$-invariance of the
Ricci tensor is the description of the  almost \ka
manifolds  $(M,g,J,\Om)$
whose sectional curvature at any point  is the same on all 
totally-real planes. We will refer to these manifolds as  almost \ka
manifolds of {\it pointwise constant 
totally-real sectional curvature}. 
Recall that a two-plane $\sigma
= X \wedge Y$ in the tangent bundle $TM$ is said to be {\it totally-real}
if $J \sigma = JX \wedge JY$ is orthogonal to $\sigma$. It was shown 
in \cite{FF,FFK} that almost \ka
manifolds of
pointwise constant totally-real sectional curvature of dimension higher
than 4 are in fact complex space forms. To our knowledge, the
4-dimensional case which we consider here is still open. 
We also note that it has been open
for some time the question of existence of a strictly almost \ka structure
of constant sectional curvature in dimension 4. Now the answer is known to
be negative, as a consequence of results of \cite{Ar3,OS2}.

The condition that $\si = X \wedge Y$ is a totally-real
plane is, in fact, equivalent to $g(JX, Y) = 0$, i.e., to $\si$ being a
Lagrangian plane with respect to the fundamental 2-form $\Om$. Thus, the
space of totally-real two-planes is determined by $\Om$ and it does not
depend
on the choice of an $\Om$-compatible almost  \ka structure.
Given a symplectic 4-manifold $(M, \Om)$, it is natural to ask if there
are any
$\Om$-compatible almost \ka structures which have pointwise constant
totally-real sectional curvature. It is easily checked that for any such a
structure the Ricci tensor is $J$-invariant 
and the Weyl tensor is self-dual (Lemma \ref{lem6} below). In the compact case  we conjecture
that {\it any compact, 4-dimensional almost \ka manifold of
pointwise constant totally-real sectional curvature is,
in fact, a self-dual \ka surface}. The compact self-dual
K{\"a}hler surfaces are completely described by \cite{B-YC} (see also \cite{Bu}, \cite{De},
\cite{It}) --- they are all Riemannian locally symmetric spaces and in
particular have constant scalar curvature. We are able to prove the
following:

\begin{theo}\label{th3}
Any compact almost \ka 4-manifold of pointwise
constant totally-real sectional curvature and of constant scalar curvature
is K{\"a}hler.
\end{theo} 

\noindent
Note also that the complete classification of compact Hermitian surfaces
of (pointwise)
constant totally-real sectional curvature was recently obtained in \cite{A,AGI1}.

\vspace{0.2cm}

The paper is organized as follows:

Section 2 is devoted to local considerations on almost
Hermitian
manifolds with $J$-invariant Ricci tensor. Many of our results here could
be considered as further
development of the so called ``{\it Riemannian Goldberg-Sachs theory}''
\cite{AG},
with applications to almost \ka 4-manifolds (see subsection 2.2).
In particular, we obtain in Corollary \ref{c1}
that {\it any almost K{\"a}hler anti-self-dual 4-manifold with $J$-invariant
Ricci tensor
is K{\"a}hler scalar-flat}.
Slightly stronger results are obtained computing the 
Bach tensor of an almost \ka 4-manifold with $J$-invariant Ricci tensor
(see Lemma \ref{lem4} in subsection 2.3 and its corollary). 
The technical part of our work uses the $U(2)$-decomposition of
the curvature, first and second Bianchi identities and some
Weitzenb{\"o}ck formulas.

In Section 3 we specify the local considerations of the previous
section to the case of almost K{\"a}hler 4-manifolds satisfying the
second Gray condition. We obtain that if such a manifold $(M, g, J, \Om)$
is {\it not} K{\"a}hler, then it admits a {\it negative} almost K{\"a}hler
structure 
$({\bar M}, g, {\bar J}, {\bar \Om})$, which has quite strong symmetry
properties of the curvature as well 
(Proposition \ref{prop2}). Imposing also compactness, the proof of our main 
result, Theorem \ref{th1},
is done in several steps: Assuming that the given almost \ka structure
is not K{\"a}hler, we first use some results from the Seiberg-Witten
theory to prove that the Euler number and the signature
of the manifold both vanish.  Then we show that the negative almost
\ka structure must be K{\"a}hler, and  using the Kodaira classification 
we conclude that the corresponding complex surface belongs to class VI. 

Section 4 is devoted to almost \ka 4-manifolds of (pointwise) constant
totally-real sectional curvature. The proof of Theorem 2 is an
application of the results of Section 2 and \cite{Dr1,Ar}.

\vspace{0.2cm}
\noindent  {\bf Acknowledgments}:   Both authors are especially
grateful to J. Armstrong  for valuable comments  and to
U. Semmelmann for pointing out a gap in the proof of Proposition
\ref{prop2}. They would also like to thank D. Kotschick for some stimulating
discussions.

\vspace{0.2cm}
\noindent {\bf Dedication}: We would like to dedicate this paper to the memory
of Alfred Gray, whom we did not know personally, but whose mathematical work
inspired us. 

\section{Almost Hermitian 4-manifolds with $J$-invariant Ricci
tensor}

\subsection{The $U(2)$-decomposition of the curvature tensor of an almost
Hermitian 4-manifold.}

Let $(M, g)$ be a 4-dimensional, oriented Riemannian manifold. The Hodge
operator $*$, acting as an involution on the bundle of 2-forms,
induces the
orthogonal splitting:
$$ \La ^2 M = \La^+ M \oplus \La^- M, $$
where $\La^+ M$, ($\La^- M$), is the bundle of self-dual, (anti-self-dual)
2-forms. We will freely identify vectors and co-vectors via the metric $g$
and,
accordingly, a tensor $\phi \in T^*M^{\otimes 2}$ with the corresponding
endomorphism of the
tangent bundle $TM$, by $\phi(X, Y) = g(\phi (X), Y)$. Also, if $\phi,
\psi
\in T^*M^{\otimes 2}$, by $ \phi \circ \psi$ we understand the
endomorphism
of $TM$ obtained by the composition of the endomorphisms corresponding to
the two tensors.

Considering the Riemannian curvature tensor $R$ as a symmetric
endomorphism of
$\La ^2 M$, we have the following well-known SO(4)-decomposition:
\begin{equation}\label{So4}
R = \frac{s}{12} id + {\tilde {Ric_0}} + W^+ + W^- ,
\end{equation}
where:
\newline
- $s$ is the scalar curvature;
\newline
- $ {\tilde {Ric_0}} $ is the Kulkarni-Nomizu extension
of the traceless Ricci tensor, $Ric_0$, to an endomorphism of $\La^2 M$
anti-commuting with $*$;
\newline
- $W^{\pm}$ are respectively the self-dual and anti-self-dual parts of
the Weyl tensor $W$.
\newline
We view $W^{\pm}$, respectively, as sections of the
bundles ${\cal W}^{\pm} = Sym_0 (\La ^{\pm}M)$ of symmetric, traceless
endomorphisms
of $\La^{\pm}M$. The manifold $(M, g)$ is said to be self-dual
(anti-self-dual), if $W^- = 0$ ($W^+ = 0$). The manifold is Einstein if
${\tilde {Ric_0}}$ (or, equivalently, $Ric_0$) vanishes identically.

\vspace{0.2cm}
\noindent
Let $(M, g, J, \Om)$ be an almost Hermitian, 4-dimensional
manifold, i.e., a 4-diminsional smooth manifold endowed with  a triple 
$(g,J,\Om )$ of a Riemannian metric $g$, an orthogonal almost complex structure $J$ and
a non-degenerate 2-form, not necessarily closed, satisfying the
compatibility relation (1).
The action of the almost
complex structure $J$ extends to the cotangent bundle $T^{*}M$, by
$$ (J\alpha) (X) = - \alpha (JX) , $$
and to the bundle of 2-forms $\La ^2 M$, by
$$ (J \phi) (X, Y) = \phi(JX, JY).$$
Because of this action, we have the following
splitting of the bundle of (real) 2-forms.
$$\Lambda^2 M = \Lambda^{inv} M + \Lambda^{anti} M =
         {\Bbb R}\Om + \Lambda_0 ^{inv} M + \Lambda^{anti} M,$$
where $\Lambda^{inv} M, \Lambda^{anti} M$ are respectively the spaces
of $J$-invariant and $J$-anti-invariant 2-forms and $\Lambda_0 ^{inv} M$
is the space of trace-free $J$-invariant 2-forms. Note that $\La^{anti} M$
is
the real underlying bundle of the anti-canonical bundle $(K_J)^{-1}:=
\La^{0,2}M$ of $(M, J)$; the induced complex structure ${\cal J}$ on
$\La^{anti} M$ is then given by
${\cal J}\phi(.,.):=-\phi(J.,.)$.
The vector bundle $\La^{inv}M$ of $J$-invariant real 2-forms  is the
real
underlying bundle of the bundle of (1,1) forms.

Whereas the
above decomposition
holds in any dimension, in dimension 4 it links very nicely with the
self-dual,
anti-self-dual decomposition of 2-forms
 given by the Hodge star operator:
\begin{equation} \label{2}
 \Lambda^+ M = {\Bbb R}\Om \oplus \Lambda^{anti} M \;,
\end{equation}
\begin{equation} \label{3}
    \Lambda^- M =  \Lambda_0 ^{inv} M .
\end{equation}

Because of the splitting (\ref{2}), the bundle ${\cal W^+}$ decomposes
into the
sum of three sub-bundles,
${\cal W}_{1}^+$, ${\cal W}_{2}^+$, ${\cal W}_{3}^+$ defined as follows,
see \cite{TV}:
\newline
- ${\cal W}_1 ^+ = {\Bbb R} \times M $ is the sub-bundle of elements
preserving the
decomposition (\ref{2}) and acting by homothety on the two factors; hence
is
the trivial line
bundle generated by the element $ \frac{1}{8} \Om \otimes \Om -
\frac{1}{12}
id^+ $.
\newline
- ${\cal W}_2 ^+ = \La^{anti} M$ is the sub-bundle of elements which
exchange the
two factors in (\ref{2}); the real isomorphism with $\La^{anti} M$ is seen
by identifying each
element $\phi$ of $\La^{anti} M$ with the element $ \frac{1}{2} (\Om
\otimes
\phi +
\phi \otimes \Om )$ of ${\cal W}_2 ^+$.
\newline
- ${\cal W}_3 ^+ = Sym_0(\La^{anti} M)$ is the sub-bundle of
elements preserving the
splitting (\ref{2}) and acting trivially on the first factor ${\Bbb R} \Om
$.

Consequently, we have the following splitting of the Riemannian curvature
operator $R$ (\cite{TV}):

\begin{equation}\label{U2}
R = \frac{s}{12} id + ({\tilde {Ric_0}})^{inv} + ({\tilde {Ric_0}})^{anti}
+
 W_1 ^+ + W_2 ^+ + W_3 ^+ + W^- ,
\end{equation}

\noindent
where $ ({\tilde {Ric_0}})^{inv} $ and $ ({\tilde {Ric_0}})^{anti} $ are
the
Kulkarni-Nomizu extensions of the $J$-invariant and the
$J$-anti-invariants
parts of
the traceless Ricci tensor, respectively, and $W_i ^+$ are the projections
of $W^+$ on
the spaces ${\cal W}_i ^+, \; i = 1,2,3$. The component $W_1^+$
 is given by
\begin{equation}\label{w^+_1}
W_1^+ = \frac{\kappa}{8} \Om \otimes \Om - \frac{\kappa}{12} id^+ ,
\end{equation}
where $\kappa $ is the so called {\it conformal scalar curvature} of $(g,
J)$,
defined by
$$ \kappa = 3 <W^+(\Om), \Om>,$$
where $<.,.>$ denotes the extension of the Riemannian metric $g$ to
the bundle $\Lambda ^2M$, so that the norm induced by
$<.,.>$ to be half of the usual tensor norm on
$\Lambda ^1M\otimes \Lambda ^1M$.
The component
$W_2^+$ is determined by the skew-symmetric part of the {\it star-Ricci
tensor}, defined for an almost Hermitian manifold by:
$$ Ric^{*} (X, Y) =  -R(\Om) (JX, Y) . $$
We thus get
\begin{equation}\label{w^+_2}
W^+_2=-\frac{1}{4}(\Psi\otimes \Omega + \Omega\otimes \Psi),
\end{equation}
where  $-\frac{1}{2}J\Psi$ is the skew symmetric part of $Ric^*$.
Note that in general, the star-Ricci tensor is neither symmetric nor
skew-symmetric. 
The {\it star-scalar curvature} $s^{*}$ is the trace of $Ric^*$.
It is easy to check the following relation between the scalar curvatures
of
a 4-dimensional almost Hermitian manifold:
\begin{equation}\label{kappa}
\kappa = \frac{1}{2} (3s^* - s) .
\end{equation}
On any almost \ka manifold (of arbitrary dimension), the difference between
the
star-scalar curvature
and the
scalar curvature ``measures'' the integrability of the almost complex
structure:
\begin{equation}\label{s-s^*}
s^{*} - s = \frac{1}{2} |\na J|^2 .
\end{equation}
Let us also remark here that
$ ({\tilde {Ric_0}})^{anti} $ determines the part of the curvature acting
from
$\La ^{anti} M$ to $\La ^{-} M$, thus a 4-dimensional
 almost Hermitian manifold has
$J$-invariant Ricci tensor if and only if 
$$ < R(\La ^{anti} M), \La ^- M > = 0.$$

\subsection{The second Bianchi identity for almost Hermitian 4-manifolds
with $J$-invariant Ricci tensor}

The co-differential $\delta W^+$ of the positive Weyl tensor of $(M, g)$
is a section of the rank 8 vector bundle

\[ {\cal V} = Ker ( tr : \La^1 M \otimes \La^+ M \longrightarrow \La^1 M )
,\]
where $ tr $ is defined by $ tr(\alpha \otimes \phi) = \phi(\alpha) $ on
decomposed elements. The vector bundle ${\cal V}$ splits as
${\cal V} = {\cal V}^+ \oplus {\cal V}^-$, see \cite{AG}, where:\\

\noindent
(a) ${\cal V}^+$
is identified with the cotangent bundle $ T^{*} M$ by

\begin{equation}
\label{alpha}
 T^{*}M \ni \alpha \mapsto A = J\alpha \otimes \Om -
\frac{1}{2}\sum_{i=1}^{4}e_{i} \otimes (\alpha \wedge e_{i} - J\alpha
\wedge
Je_{i}) ,
\end{equation}

\[ {\cal V}^{+} \ni A \mapsto \alpha = - \frac{1}{2}J<A,\Om> ,\]

\noindent where $<A,\Om>$ denotes the 1-form defined by $
X \mapsto <A_{X},\Om >$. \\

\noindent (b) $\cal V^{-}$ is identified, as a real vector bundle,  to the
real, rank 4,  vector bundle  $\Lambda^{0,1}M \otimes K^{-1}_{J}$; if $\phi$ is a
non-vanishing local
section of $\Lambda^{anti} M$, then it trivializes $K_J^{-1}$ and ${\cal
V}^-$
can be again identified with $T^*M$  by
\begin{eqnarray}\label{beta}
\beta \in T^{*}M &\mapsto &B = \sum_{i=1}^{4}e_{i} \otimes (J\beta \wedge
\phi(e_{i}) + \beta \wedge {\cal J}\phi(e_{i}))\\ \nonumber
B \in \cal V ^{-} &\mapsto & \beta =- \frac{1}{\vert \phi \vert
^{2}}J<B,\phi>. \nonumber
\end{eqnarray}

\vspace{0.2cm}
\noindent
We denote by  $(\delta W^{+})^{+}$ , resp. $(\delta W^{+})^{-}$, the
component of $\delta W^{+}$ on  $\cal V^{+}$, resp. on  $\cal V^{-}$.
The {\it Cotton-York tensor} $C$ of $(M, g)$ is  defined by:
$$ C_{X,Y,Z} = \frac{1}{2} \Big[ \na_{Z} (\frac{s}{12} g + Ric_0) (Y,X) -
\na_{Y} (\frac{s}{12} g + Ric_0)(Z,X) \Big] .$$
The second Bianchi identity reads as
\begin{equation} \label{bianchiII}
\delta W = C,
\end{equation}
where $\delta W$ is the
co-differential of the Weyl tensor $W$. In particular, we have
\begin{equation}\label{half}
\delta W^+ = C^+,
\end{equation}
\noindent
where $C^+$ denotes the self-dual part of $C_X$, $X \in TM$.

\vspace{0.2cm}

\noindent
For an almost Hermitian 4-manifold $(M, g, J, \Om)$, denote by $\theta$
the
Lee 1-form
given by $\theta = J \delta \Om$. Then we have the following:

\vspace{0.2cm}

\noindent
\begin{Lemma}\label{lem1}
Let $(M, g, J, \Om)$ be an almost Hermitian 4-manifold
with
$J$-invariant Ricci tensor. Denote by $\alpha$ the 1-form corresponding
to $ (\delta W^+)^+ $ via the isomorphism (\ref{alpha}).
Then $\alpha$ is given by
$$\alpha = -\frac{ds}{12} + \frac{1}{2} Ric_0 (\theta).$$
\end{Lemma}
{\it Proof:} From (\ref{half}) we have
$$ \alpha(X) = - \frac{1}{4} \sum_{i=1}^{4} \na_{e_{i}} (\frac{s}{12} g +
Ric_0)(Je_i, JX) = $$
$$ = - \frac{1}{4} \Big[ \frac{ds}{12}(X) - (\delta Ric_0)(X) +
\sum_{i=1}^{4} Ric_0(e_i, J(\na_{e_i} J)(X)) - Ric_0(\theta, X) \Big] =$$
$$ = - \frac{1}{4} \Big[ \frac{ds}{3}(X) +
\sum_{i=1}^{4} Ric_0(e_i, J(\na_{e_i} J)(X)) - Ric_0(\theta, X) \Big] .$$
For any almost Hermitian manifold $\na J$ is given by (cf. e.g.
\cite{k-n}):
$$ \na _X J = \frac{1}{2}(X\wedge J\theta + JX\wedge \theta) +\frac{1}{2}
N_{JX}, $$ where $N_{X}(.,.)=g([J.,J.]-[.,.]-J[J.,.]-J[.,J.],X)$ is the
Nijenhuis tensor of $(M,J)$.
Using this formula and the fact that the Ricci tensor is $J$-invariant, we
compute
$$ \sum_{i=1}^{4} Ric_0(e_i, J(\na_{e_i} J)(X)) = - Ric_0(\theta, X) ,$$
and we obtain the expression claimed for $\alpha$. $\square$

\vspace{0.2cm}

\noindent
Regarding the component $(\delta W^+)^- $, let $\phi$ be a (locally
defined) section of
$\Lambda ^{anti}M$, with $|\phi|^2=2$, and denote by $\beta$  the
1-form corresponding to $(\delta W^+)^-$ via the isomorphism
(\ref{beta}). 
It follows from (\ref{alpha}) and (\ref{beta}) that 
$$\beta =\frac{1}{2}\Big( -J<\delta W^+,\phi> - \frac{1}{2}{\cal
J}\phi(\alpha) \Big) .$$
To compute $J<\delta W^+, \phi>$ we will proceed in the same way as
computing $J<\delta W^+, \Om>$ in the proof of Lemma \ref{lem1}; we
will now consider instead of $J$ the almost complex structure $I$ 
whose K{\"a}hler form is $\phi$. Observe that $Ric_0$ is now
$I$-anti-invariant. Involving Lemma \ref{lem1}, together with the
folowing 
expressions for the covariant derivatives of $\na J$ and $\na I$:
$$\na J = a \otimes I + b \otimes (J\circ I);$$
$$\na I = -a\otimes J + c\otimes (J\circ I), $$
where $a, b, c$ are some 1-forms, we eventually obtain 
\begin{equation}\label{beta0}
\beta = -\frac{1}{4}Ric_0(a + Jb).
\end{equation}  
Note that the almost complex structure $J$ is integrable (resp. almost
K{\"a}hler) iff $b=Ja$ (resp. $b=-Ja$); we thus obtained the following:
\begin{Lemma}\label{lem2} {\rm {\bf (\cite[Proposition 4]{AG})}}
Let $(M, g, J, \Om)$ be an almost Hermitian
4-manifold with $J$-invariant Ricci tensor. Suppose that the traceless
Ricci
tensor of $g$ does not vanish. Then $J$ is integrable if and only if
$$ (\delta W^+)^- = 0 .$$
\end{Lemma}

\begin{cor}\label{c1}
Any anti-self-dual almost \ka 4-manifold of $J$-invariant
Ricci tensor is K{\"a}hler (hence, scalar flat).
\end{cor}
{\it Proof:} Suppose that the traceless Ricci tensor does not vanish on an
open
subset $U$ of $M$. It follows from Lemma \ref{lem2} that $(g, J, \Om)$ is
K{\"a}hler on $U$.
Hence
the scalar curvature vanishes on $U$ and, by continuity, this holds on the
closure
of $U$. If $M - {\bar U}$ is
 non-empty, then $(g, J, \Om)$ is an Einstein, anti-self dual, almost \ka
structure on
the open set
$M - {\bar U}$. By \cite[Theorem 2.4]{Ar3}, or by \cite[Theorem 1]{OS2},
it
follows that
$(g, J, \Om)$ is \ka on $M - {\bar U}$. Thus the scalar curvature is zero
everywhere on $M$ (see (\ref{w^+_1}), (\ref{kappa}) and (\ref{s-s^*})),
hence $(g, J, \Om)$ is a \ka structure on $M$. $\square$ 
 
\subsection{The Bach tensor of almost Hermitian 4-manifolds with
$J$-invariant Ricci tensor}

The Bach tensor $B$ of a (4-dimensional) Riemannian manifold $(M,g)$ is
defined by
\begin{equation}\label{Bach}
 B_{X,Y}=\sum_{i=1}^4 [\na_{e_i}(\delta W)(X,e_i,Y) + W(X,e_i,h(e_i),Y)],
\end{equation}
where $W$ is the Weyl conformal tensor and 
$h=\frac{1}{2}(Ric-\frac{s}{6}g)$ is
the normalized Ricci tensor of $(M,g)$.
It is known that on a compact 4-manifold
$B_{X,Y}$ is the gradient of the Riemannian
functional $g\mapsto \int_M |W|^2 dV_g$,
acting on the space of all Riemannian metrics on $M$,
cf. \cite{bach}. It thus follows that $B$
is symmetric, traceless and conformally invariant
(1,1)-tensor
on $M$, cf. \cite{bach}; it vanishes if the metric $g$ is (locally) conformal
to
Einstein by means of the second Bianchi identity (\ref{bianchiII}), or if
$g$
is self-dual or anti-self-dual by means  of the following expression for
$B$,
obtained by Gauduchon \cite{Ga}:
\begin{equation}\label{gauduchon}
\frac{1}{2}B_{X,Y}=\sum_{i=1}^4 [\na_{e_i}(\delta
W^{\pm})(X,e_i,Y)+W^{\pm}(X,e_i,h(e_i),Y)].
\end{equation}
In fact, identifying freely the $(2,0)$-tensors with the corresponding
1-forms
with values in $TM$ via the metric $g$
and using the second
Bianchi identity (\ref{bianchiII}), we get from (\ref{Bach})
\begin{equation}\label{*}
(*B)_{X,Y,Z}=d^{\na}(*C)_{X,Y,Z} + \sigma_{X,Y,Z}[(*\circ W)_{X,Y,h(Z)}],
\end{equation}
where $*:\Lambda^k M\otimes TM \mapsto \Lambda^{4-k}M\otimes TM$
and $d^{\na}:\Lambda^k M\otimes TM \mapsto \Lambda^{k+1}M\otimes TM$ are
the Hodge operator and the Riemannian differential, respectively,
acting on the bundle of the $k$-forms with values in $TM$;
$C=-d^{\na}h$ is the the Cotton-York tensor, viewed as a section of
$\Lambda^2M\otimes TM$;
$*\circ W = W^+ - W^-$ is  considered as a section of
$(\Lambda ^2 M\otimes \Lambda^1 M)\otimes TM$ and
$\sigma_{X,Y,Z}$ denotes the cyclic sum on $X,Y,Z$.\\
On the other hand, on any Riemannian manifold we have
$$ d^{\na}C_{X,Y,Z} =-(d^{\na})^2 h_{X,Y,Z}
=-\sigma_{X,Y,Z}[R_{X,Y,h(Z)}],$$
hence, using (\ref{So4}) we infer
$$d^{\na}C_{X,Y,Z} + \sigma_{X,Y,Z}[W_{X,Y,h(Z)}]=0, $$
which, together with (\ref{*}), implies
\begin{eqnarray} \label{**}
(*B)_{X,Y,Z}&=& 2({d^{\na}}(C^+)_{X,Y,Z} +
\sigma_{X,Y,Z}[W^+_{X,Y,h(Z)}])\\ \nonumber
 & =& -2(d^{\na}(C^-)_{X,Y,Z} + \sigma_{X,Y,Z}[W^-_{X,Y,h(Z)}]).\nonumber
\end{eqnarray}
Relation (\ref{gauduchon}) follows now by applying the Hodge operator $*$ 
to both sides in (\ref{**}).

\vspace{0.2cm}

The Bach tensor can be also expressed in terms of the Ricci tensor $Ric$
as
follows:
Using the second Bianchi identity (\ref{bianchiII}), relation (\ref{Bach})
can be
rewritten as
\begin{equation}\label{intelect}
B = \delta^{\na}d^{\na} h - \stackrel{\circ}{W}(h),
\end{equation}
where:
\newline
- $h$ is considered as a section of $\La ^1 M\otimes TM$;
\newline
- $\delta^{\na}:\Lambda^{k+1} M\otimes TM \mapsto \Lambda^{k}M\otimes TM$
is the
formal adjoint of $d^{\na}$;
\newline
- $\stackrel{\circ}{W}(t)= -\sum_{i=1}^4W(X,e_i,t(e_i),Y)$ is the action
of $W$
on
$\La ^1 M\otimes TM$.\\
On the other hand, the Weitzenb{\"o}ck formula for the 1-forms with values
in
$TM$
reads as \cite[Proposition 4.1]{Bu}
$$(\delta^{\na}d^{\na} + d^{\na}\delta^{\na})t = \na^*\na t + t\circ Ric -
\stackrel{\circ}{R}(t),$$
where $t$ is a section of $\La ^1 M\otimes TM$ and
$\stackrel{\circ}{R}(t)=
-\sum_{i=1}^4R(X,e_i,t(e_i),Y)$. Using the above formula, specified for
$h=\frac{s}{24}g + \frac{1}{2}Ric_0$, together with the Ricci identity
$\delta^{\na}h=-\frac{ds}{6}$, we get from (\ref{intelect}):
\begin{equation}\label{bour}
B = \na^*\na h  +\frac{1}{6}\na ds + h\circ Ric
-\stackrel{\circ}{W}(h) -\stackrel{\circ}{R}(h).
\end{equation}
If $(M,g,J)$ is an almost Hermitian 4-manifold with $J$-invariant Ricci
tensor, then the traceless Ricci tensor has eigenvalues
$(\lambda, \lambda, -\lambda, -\lambda)$ and hence
$$h\circ Ric = \Big( \frac{s^2}{96} + \frac{|Ric_0|^2}{8} \Big) g + 
\frac{s}{6}Ric_0,$$
and by (\ref{So4})
$$\stackrel{\circ}{R}(h)= \Big( \frac{s^2}{96}+\frac{|Ric_0|^2}{8} \Big) g
+ 
\stackrel{\circ}{W}(h).$$
Thus,  (\ref{bour}) reduces to (see also \cite[(24) and Lemma 4,(i)]{De}):
\begin{equation} \label{intelect2}
B = \frac{1}{2}\na^*\na Ric_0  + \frac{\Delta s}{24}g +
\frac{1}{6}\na ds + \frac{s}{6}Ric_0 - \stackrel{\circ}{W}(Ric_0).
\end{equation}
If moreover $(M,g,J)$ is an {\it almost K{\"a}hler} 4-manifold, we use
(\ref{intelect2}) to obtain the following formula for $B$, which has a
close
resemblance with the expression of the Bach tensor on a 4-dimensional
K{\"a}hler manifold obtained in \cite{De}.

\begin{Lemma}\label{lem4}
Let $(M,g,J)$ be an almost \ka 4-manifold with $J$-invariant Ricci tensor.
Then the Bach tensor $B$ of $g$ is given by
$$ B = -\frac{1}{3}(\na ds)^{inv} +\frac{1}{6}(\na ds)^{anti} -
\frac{\Delta
s}{12}g - \frac{s}{6}Ric_0 + S,$$
where:

- $(\na ds)^{inv}$ and $(\na ds)^{anti}$ denote the $J$-invariant and
$J$-anti-invariant parts of the symmetric $(2,0)$-tensor field $\na ds$,
respectively;

- $S$ is the symmetric, $J$-anti-invariant tensor, given by:
$$ S = \sum_{i=1}^{4}(\na_{e_i} {\rho}_0) \circ (\na_{e_i} \Om),$$
with $\rho_0(.,.)=Ric_0(J.,.)$ being the traceless Ricci form.
\end{Lemma}
{\it Proof:} Since $Ric_0$ is $J$-invariant, we have
\begin{eqnarray} \nonumber
\na^*\na Ric_0 &=& -\na^*\na (\rho_0 \circ \Om) \\ \label{tedi1}
               &=& -(\na^*\na \rho_0)\circ \Om -\rho_0 \circ (\na^*\na
\Om) + 2S.\nonumber
\end{eqnarray}
To compute $\na^*\na \Om$ and $\na^*\na \rho_0$ we will use the
Weitzenb{\"o}ck
formula for 2-forms (also considered as sections of
$\Lambda^1 M\otimes TM$); for any 2-form $\phi$ we have (cf. e.g.
\cite{Be}):
$$\Delta \phi = \na^*\na \phi +\frac{s}{3}\phi -2W(\phi). $$
Since the K{\"a}hler form $\Om$ is harmonic, we get
\begin{equation} \label{tedi2}
\rho_0\circ (\na^*\na \Om)= \frac{s}{3}Ric_0  + 2\rho_0\circ W^+(\Om).
\end{equation}
Moreover, it is known that for an almost K{\"a}hler 4-manifold with
$J$-invariant Ricci tensor,
the Ricci form $\rho = \frac{s}{4}\Om + \rho_0$ is closed, see \cite{Dr2}.
We thus compute
$$\delta \rho_0 = J(\delta Ric_0)= -\frac{Jds}{4}, $$
$$d\rho_0 =-\frac{ds}{4}\wedge \Om,$$
$$\Delta \rho_0 = -(\na ds)^{inv}\circ \Om -\frac{\Delta s}{4}\Omega,$$
hence
\begin{equation} \label{tedi3}
(\na^*\na \rho_0) \circ \Om =  (\na ds)^{inv}  + \frac{\Delta s}{4}g +
\frac{s}{3}Ric_0 + 2W^-(\rho_0)\circ \Om.
\end{equation}
Since $h=\frac{s}{24}g + \frac{1}{2}Ric_0$ is $J$-invariant, we
have
$$ \sum_{i=1}^4 W^{+}_3(X,e_i,h(e_i),Y)=0. $$
Thus, using (\ref{w^+_1}) and
(\ref{w^+_2}) we compute
\begin{eqnarray}\nonumber
\stackrel{\circ}{W^{+}}(h) & = &
\frac{1}{2}\stackrel{\circ}{W^{+}}(Ric_0) \\ \nonumber
 &=&
\frac{\kappa}{12}Ric_0 -
\frac{1}{8}[{\cal J}\Psi\circ Ric_0-Ric_0\circ {\cal J}\Psi].
\end{eqnarray}
As $Ric_0$ is $J$-invariant, it anti-commutes with 
${\cal J}\Psi$ and the latter
equality reduces to
\begin{equation}\label{tedi0}
\stackrel{\circ}{W^{+}}(Ric_0)= \frac{\kappa}{6}Ric_0 -
\frac{1}{2}({\cal J}\Psi)\circ Ric_0.\nonumber
\end{equation}
Using now (\ref{w^+_1}), (\ref{w^+_2}) and (\ref{tedi0}), we infer
$$2\rho_0\circ W^+(\Om)=\rho_0\circ(\frac{\kappa}{3}\Om -\Psi)=
-2\stackrel{\circ}{W^+}(Ric_0).$$
Similarly, considering instead of $J$
the {\it negative} almost Hermitian structure $\bar{J}$ that makes the
Ricci tensor ${\bar J}$-invariant (at the points where $Ric_0\neq 0$), we get
$$2W^{-}(\rho_0)\circ\Om = - 2\stackrel{\circ}{W^-}(Ric_0).$$
Hence, substituting (\ref{tedi2}) and (\ref{tedi3}) into (\ref{tedi1}),
we eventually obtain
$$\frac{1}{2}\na^*\na Ric_0 = S -\frac{1}{2}(\na ds)^{inv} -\frac{\Delta
s}{8}g
-\frac{s}{3}Ric_0  + \stackrel{\circ}{W}(Ric_0) ,$$
and the claim follows by (\ref{intelect2}). $\square$

\begin{cor}\label{c2}
Let $(M,g,J,\Omega)$ be an almost K{\"a}hler 4-manifold with $J$-invariant
Ricci
tensor and constant scalar curvature. Suppose that the Bach tensor
vanishes.
Then either $g$ is Einstein, or else the scalar curvature of $g$ is zero.
\end{cor}
{\it Proof:} Considering the $J$-invariant part of $B$ given by Lemma
\ref{lem4},
we get $sRic_0=0$ and the claim follows. $\square$

\vspace{0.2cm}
\noindent
{\bf Remark 1.} Corollary \ref{c2} can be considered as a generalization
of
Corollary \ref{c1}. Indeed, if $(M,g,J)$ is an anti-self-dual
almost-K{\"a}hler
4-manifold with $J$-invariant Ricci tensor, then the scalar curvature is
constant by Lemma \ref{lem1} and the Bach tensor vanishes by
(\ref{gauduchon}); moreover, it follows from (\ref{w^+_1}),
(\ref{kappa}) and (\ref{s-s^*})
that the scalar curvature of an anti-self-dual almost \ka metric $g$ is
zero if and only if $(g,J)$ is K{\"a}hler. $\Diamond$

\vspace{0.2cm}
\noindent
{\bf Remark 2.} 
Using Lemma \ref{lem4} in the compact case, a stronger result than
Corollary \ref{c2} could be obtained. In fact {\it any 
compact almost K{\"a}hler 4-manifold with
$J$-invariant Ricci
tensor and non-positive scalar curvature, whose Bach tensor
vanishes is either an Einstein manifold, or else it is a K{\"a}hler,
scalar-flat surface.}  Indeed, taking the inner product of the
Bach tensor, as given in 
Lemma \ref{lem4}, with $Ric_0$ and integrating over the manifold, we get:

$$ 0 = \int_M \Big( - \frac{1}{3} <\na d s, Ric_0> - \frac{s}{6} |Ric_0|^2
\Big) dV_g .$$
Integrating the first term by parts and taking into consideration that
$\delta Ric_0 = \frac{1}{4} ds$, the above relation changes to:
$$ 0 = \int_M \Big( \frac{|ds|^2}{12}  - \frac{s}{6} |Ric_0|^2 \Big) dV_g
.$$
If $g$ is {\it not} Einstein the above formula shows that
$s$ must identically vanish 
and then the conclusion follows from \cite[Theorem 1]{Dr1}. $\Diamond$

\vspace{0.2cm}
\noindent
Note that both Corollaries \ref{c2} and Remark 2 apply 
for almost K{\"a}hler, self-dual manifolds with $J$-invariant Ricci
tensor,
as well. 

\vspace{0.2cm}
\noindent
{\bf Remark 3.}
It is well known that compact quotients of the complex hyperbolic space
${\bf C}{\cal H}_2$ admit K{\"a}hler-Einstein, self-dual metrics of 
negative scalar curvature. In the 
conformal class of such a metric strictly almost \ka metrics could be
found, as it follows from \cite{ApDr}. 
Thus, there are examples of self-dual, conformally Einstein, strictly
almost 
\ka metrics on compact 4-manifolds. On the other hand, a simple Bochner type 
argument shows
that there are no strictly almost \ka metrics in the conformal class
of an anti-self-dual (equivalently, scalar-flat), \ka metric 
(see {\it e.g.} \cite{ApDr}). However, it was pointed out  by J. Armstrong
\cite{Ar} that strictly almost \ka anti-self-dual metrics do exist on  
generic ruled surfaces. This follows from the analysis of the moduli
spaces of anti-self-dual metrics and of scalar-flat, \ka metrics done
by King, Kotschick \cite{KK} (see also \cite{It}) and
LeBrun, Singer \cite{Le-Si}, respectively. $\Diamond$

\section{The second curvature condition of Gray in dimension 4}

\vspace{0.2cm}

In \cite{Gr}, A. Gray considered almost Hermitian manifolds whose
curvature
tensor has a certain degree of resemblance to that of a \ka manifold.
The following identities arise naturally:
\newline $ G_1 \; \; \; \; R_{XYZW} = R_{XYJZJW} $
 ;
\newline $ G_2 \; \; \; \; R_{XYZW} - R_{JXJYZW} = R_{JXYJZW} + R_{JXYZJW}
$
 ;
\newline $ G_3 \; \; \; \; R_{XYZW} = R_{JXJYJZJW} $
 .
\newline
We will call the identity $G_i$ as the i-th condition on the curvature of
Gray.
It is a simple application of the first Bianchi identity
to see that $G_1 \Rightarrow G_2 \Rightarrow G_3 $. Also elementary is
the fact that a K{\"a}hler structure satisfies relation $G_1$ (hence, all
of the
relations $G_i$). Following \cite{Gr}, if ${\cal AK}$ is the class of
almost
\ka
manifolds, let ${\cal AK}_i$ be the subclass of manifolds whose curvature
satisfies
identity $G_i$. We have the obvious inclusions
$$ {\cal AK} \supseteq {\cal AK}_{3} \supseteq {\cal AK}_2 \supseteq
   {\cal AK}_1 \supseteq {\cal K} , $$
where ${\cal K}$ denotes the class of \ka manifolds. In \cite{Gr}
it was proven that the equality $ {\cal AK}_1 = {\cal K}$ holds locally
(see also \cite{Go}),
and it was also shown that the inclusion $ {\cal AK} \supset {\cal AK}_3$
is strict.   The examples of Davidov
and Mu\u{s}karov \cite{DM}, multiplied by compact \ka manifolds, show
that
even in the compact case,  the inclusion $ {\cal AK}_2 \supset {\cal K}$
is strict in dimension $2n \ge 6$.

In the compact, 4-dimensional case, it becomes apparent that the topology
of the underlying manifold
has consequences on the relationships between these classes. It was shown
in
\cite[Theorem 3]{Dr2} that for a compact 4-manifold with second Betti
number
equal to 1, the equality ${\cal AK}_3 = {\cal K}$ holds.
We will now deal with almost \ka 4-dimensional manifolds satisfying the
second curvature condition of Gray.

First, let us observe that the condition $G_2$ can be expressed in terms of 
the  $U(2)$-decomposition (\ref{U2}) of the curvature as follows:

\begin{Lemma}\label{lem5}
An almost Hermitian 4-manifold $(M, g, J, \Om)$ satisfies
the second curvature condition of Gray
if and only if the Ricci tensor is $J$-invariant, $W_2 ^+ = 0$ and
$W_3 ^+ = 0$.
\end{Lemma}
{\it Proof:} Easy consequence of (\ref{U2}), see \cite{TV}. $\square$

\begin{prop}\label{prop1}
For an almost Hermitian 4-manifold satisfying the condition $G_2$, we have
$$d(s - s^{*}) - \kappa \theta = 4 Ric_0 (\theta).$$
\end{prop}
{\it Proof:} Using Lemma \ref{lem5}, the positive Weyl tensor of an
almost Hermitian
4-manifold of the class ${\cal AH}_2$ is given by (\ref{w^+_1}).
Thus, computing directly,
$$ \alpha = - \frac{1}{2} J<\delta W^+, \Om> = -\frac{1}{8} \kappa \theta
- \frac{1}{12} d\kappa .$$
From this, using Lemma \ref{lem1} and (\ref{s-s^*}), we obtain the
identity
claimed. $\square$

\vspace{0.2cm}

\noindent
{\bf Remark 4.} It follows from the Riemannian version of the
Robinson-Shild
theorem in General Relativity that for the class of Hermitian 4-manifolds,
the
condition that $Ric$ is $J$-invariant is in fact equivalent to the
condition $G_2$, \cite[Theorem 2]{AG}. Thus, Proposition
\ref{prop1}
holds for
any Hermitian surface of $J$-invariant Ricci tensor. $\Diamond$

\vspace{0.2cm}

\noindent
Since on an almost \ka 4-manifold, the 1-form $\theta$ vanishes, we have
the
following immediate consequence  of Proposition \ref{prop1}.

\begin{cor}\label{c3}
 Let $(M, g, J, \Om)$ be a 4-dimensional almost \ka manifold that
satisfies the condition $G_2$. Then $ |\na J|^2 = 2( s^{*} - s)$ is a
constant.
\end{cor}
It follows from Corollary \ref{c3} and (\ref{s-s^*}), that if
$(M,g,J,\Omega)$
is almost K{\"a}hler, {\it non-K{\"a}hler} 4-manifold in the class ${\cal
AK}_2$,
then the K{\"a}hler nullity
$${\cal D} = \{ X \in TM| \na_X J = 0 \} $$
of $(g,J)$ is a well-defined 2-dimensional distribution over $M$. If we
denote
${\bar M}$ the manifold $M$ with the reversed orientation, then
we may consider the $g$-orthogonal almost complex structure ${\bar J}$ on
${\bar M}$, defined in the following manner: ${\bar J}$ coincides with $J$
on
${\cal D}$ and  ${\bar J}$ is equal to $-J$ on
${\cal D}^{\perp}$. Denote by ${\bar \Om}$ the fundamental form of
$(g, {\bar J})$ and by $W^-_i, i=1,2,3$ the $U(2)$ components of the
negative Weyl tensor $W^-$ determined by $(g,{\bar J})$.
Then we have the following:

\begin{prop}\label{prop2}
Let $(M, g, J, \Om)$ be a 4-dimensional almost K{\"a}hler, {\rm
non-K{\"a}hler}
manifold, satisfying the condition $G_2$. Then
\begin{enumerate}
\item[{\rm (i)}] The traceless Ricci tensor $Ric_0$ of $g$ is given by
$$ Ric_0=\frac{\kappa}{4}[-g^{\cal D} + g^{{{\cal D}}^{\perp}}],$$
where $g^{\cal D}$ (resp. $g^{{\cal D}^{\perp}}$) denotes the restriction
of $g$ on  ${\cal D}$ (resp. on ${\cal D}^{\perp}$).
\item[{\rm (ii)}] The triple $(g,{\bar J},{\bar \Om})$ is an almost
K{\"a}hler structure on ${\bar M}$  with ${\bar J}$-invariant Ricci tensor
and $W_2^- = 0$. Moreover, ${\cal D}^{\perp}$ belongs to the
K{\"a}hler nullity of $(g,{\bar J}, {\bar \Om})$. 
\end{enumerate}
\end{prop}
{\it Proof:} (i) 
With the notations of Section 2.3, let $\phi$ be a non-vanishing (local)
section of
$\Lambda ^{anti} M$, which satisfies $|\phi|^2=2$ at
any point. As we have already mentioned, for an almost \ka 4-manifold,
the covariant derivative $\na \Om$ of the K{\"a}hler form
$\Om$ can be written as
\begin{equation}\label{AK}
\na \Om =a\otimes \phi - Ja\otimes {\cal J}\phi,
\end{equation}
where the 1-form $a$ is of constant length $\frac{s^*-s}{4}$
(see (\ref{s-s^*}) and Corollary 3) and $\{a, Ja\}$ span
${\cal D}^{\perp}$. It follows from Lemma \ref{lem2} that $\beta =
-\frac{1}{2}Ric_0(a)$. Moreover,
if the manifold also satisfies the condition $G_2$, i.e., 
$W^+_2=0, W^+_3=0$ (see Lemma \ref{lem5}),  then by using (\ref{w^+_1}) and
(\ref{AK})
we directly compute:
$$\beta =\frac{1}{2}(-J<\delta W^+,\phi> + \frac{1}{2}\phi
<\delta W^+,\Om>)=-\frac{\kappa}{8}a .$$ Comparing the two expressions
for $\beta$ we get  the first part of Proposition \ref{prop2}. 

\vspace{0.2cm}

\noindent
(ii) From the Ricci identity we get
\begin{equation}\label{ric-id}
(\na^2_{X,Y}-\na^2_{Y,X})(\Om)(Z,T)= -R_{XYJZT}-R_{XYZJT}.
\end{equation}
On the other hand, it follows from (\ref{AK}) that
\begin{equation}\label{ric_id}
\na^2|_{\Lambda^2M}\Om =(da -Ja\wedge \xi)\otimes\phi - (d(Ja) + a\wedge
\xi)\otimes {\cal J}\phi,
\end{equation}
where the 1-form $\xi$ is defined from the equality $\na \phi= -a\otimes
\Om
+ \xi \otimes {\cal J}\phi$.
We thus obtain from (\ref{ric-id}) and (\ref{ric_id})
$$ da - Ja\wedge\xi = - R({\cal J}\phi); \ d(Ja) + a\wedge \xi = - R(\phi).$$
As our manifold satisfies the condition $G_2$, we easily deduce from 
Lemma \ref{lem5}, (\ref{U2}) and (\ref{kappa}) that 
$R(\phi)= \frac{(s-s^*)}{8}\phi$ and  
$R({\cal J}\phi)=\frac{(s-s^*)}{8}{\cal J}\phi$, hence the above 
equalities reduce to
\begin{equation}\label{(b)}
da = Ja\wedge \xi + \frac{(s^*-s)}{8}{\cal J}\phi; \ d(Ja)= - a\wedge \xi +
\frac{(s^*-s)}{8}\phi.
\end{equation}
We thus get from (\ref{(b)}) that $d(a\wedge Ja)=0$ and  using the fact
that
$|a|^2=\frac{(s^*-s)}{4}=const$, we conclude that the K{\"a}hler form ${\bar
\Om}=\Om -
\frac{2}{|a|^2}a\wedge Ja$ is closed, i.e., $(g,{\bar J})$ is almost
K{\"a}hler. \\
The Ricci tensor
is ${\bar J}$-invariant because of Proposition \ref{prop2}, (i). \\
We shall further use the following implications, which are immediate 
consequences of (\ref{ric-id}):

\vspace{0.2cm}
\begin{enumerate}
\item[(a)] $\na^2|_{\Lambda^-M} \Om =0$ if and only if the Ricci tenor is $J$-invariant;

\vspace{0.2cm}
\item[(b)] $\na^2|_{\Om} \Om =0 $ if and only if $W^+_2=0$;

\vspace{0.2cm}
\item[(c)] $(\na^2_{Z_1,Z_2}-\na^2_{Z_2,Z_1})(\Om)(Z_3,Z_4)=0, \forall
 Z_1,Z_2,Z_3,Z_4\in T^{1,0}M$ if and only if $W^+_3=0$.

\end{enumerate}

\vspace{0.1cm}
\noindent
To see that
$W^-_2=0$, we will prove first that $(\na_X {\bar J})=0$ for any vector
$X$ from ${\cal D}^{\perp}$. Indeed, put
$Z_1=A-iJA, Z_{\bar 1}=A+iJA; \ Z_2=B-iJB, Z_{\bar 2}=B+iJB$, where $\{A,
JA\}$ is an orthonormal frame of ${\cal D}^{\perp}$ so  as $A$ and $JA$ to
be the dual orthonormal frame of $\{ a,Ja\}$, and $\{ B,JB\}$ is an
orthonormal frame of
${\cal D}$. Since $(g,{\bar J})$ is almost K{\"a}hler, it is sufficient to prove
that $(\na_{Z_1}{\bar \Om})(Z_1,Z_{\bar 2})=0$. Because
$(\na_{Z_1}\Om)(Z_1,Z_{\bar 2})=0$ this is equivalent to
$\na_{Z_1} (a\wedge Ja)(Z_1,Z_{\bar 2})=0$;
the latter equality follows from (\ref{(b)}). Now,
using that the Ricci tensor is also ${\bar J}$-invariant 
and the fact that $\{ A,JA \}$ is involutive, we easily
conclude
$$\na^2|_{\bar \Om} {\bar \Om}=\na^2|_{\Om} {\bar \Om} + 2(\na^2_{A,JA}
-\na^2_{JA,A}){\bar \Om}=0 .$$  
Thus we obtain $W^-_2=0$ by the corresponding 
version of the equivalence (b) written with respect to the negative almost
K{\"a}hler 
structure  ${\bar \Omega}$.  $\square$

\vspace{0.1cm}

We are ready now to prove Theorem 1.

\vspace{0.1cm}
\noindent
{\bf Proof of Theorem 1:} Let us assume that $(M, g, J, \Om)$ is a compact, 
strictly almost \ka 4-dimensional manifold, satisfying the second curvature condition of Gray. 
Then by Corollary \ref{c3} we know that $ s^{*} - s $ is a non-zero constant function. 
Armstrong \cite{Ar2} showed
that a compact
4-dimensional almost \ka manifold with $ s^{*} - s$ nowhere vanishing,
satisfies the topological condition:
\begin{equation}\label{temp0}
2c_1^2(M) + c_2(M) = 0.
\end{equation}
Equivalently we have,
\begin{equation} \label{temp00}
5\chi(M) + 6 \sigma(M) = 0,
\end{equation}
since $ c_1^2(M) = 2 \chi(M) + 3\sigma(M)$ and $c_2(M) =
\chi(M)$, where $\chi(M)$ and $\sigma(M)$ are the Euler number and the
signature of $M$, respectively. \\We first prove that
the Euler number and the signature are both zero under our assumptions.
By Proposition \ref{prop2}, the manifold with the reversed orientation
${\bar M}$ also admits a symplectic form ${\bar \Om}$. 
We use some consequences of the Seiberg-Witten theory. There are
several cases, depending on the values of the Betti
numbers $b^+(M)$, $b^-(M)$, where $b^+(M)$ ($b^-(M)$) is equal to the
dimension of the
space of
harmonic (anti-)self-dual 2-forms. It is a classical result that $b^+(M),
b^-(M)$
do not depend on the metric and they are in fact topological invariants
of the compact 4-manifold.

\vspace{0.2cm}

\noindent
(a) $ b^+(M) > 1, \; \; b^-(M) > 1$
\newline
In this case we first remark that $M$ and ${\bar M}$ are minimal
symplectic
manifolds. Indeed, as shown in \cite{Ko}, for example, it follows that
there
are  no embedded spheres of self-intersection $\pm 1$ because the
Seiberg-Witten invariants are
non-vanishing on both $M$ and ${\bar M}$, \cite{Ta2}. By a result of
Taubes
\cite{Ta1}, a minimal compact symplectic 4-manifold $(M, \Om)$ with
$b^+(M) > 1$
satisfies
$$c_1 ^2 (M) \geq 0 .$$
Thus, we have $c_1^2(M) \geq 0$ and $ {\bar c}_1^2({\bar M}) \geq 0$,
where ${\bar c}_1$ denotes the first Chern class of $({\bar M}, {\bar
J})$.
In terms of the Euler class and the signature
of $M$:
\begin{equation}\label{temp1}
2\chi(M) + 3\sigma(M) \geq 0, \; \; 2\chi(M) - 3\sigma(M) \geq 0.
\end{equation}
Relations (\ref{temp00}) and (\ref{temp1}) imply $\chi(M) = \sigma(M) =
0$,
which is equivalent to $c_1^2(M) = c_2(M) = 0$.

\vspace{0.2cm}

\noindent
(b) $b^+(M) = 1$, $b^-(M) = 1$
\newline
Then, trivially, $\sigma(M) = b^+(M) - b^-(M) = 0$ and
we conclude from (\ref{temp00}) that $c_1^2(M) = c_2(M) = 0$ again.

\vspace{0.2cm}

\noindent
(c) $b^+(M) = 1$, $b^-(M) > 1$
\newline
We will show that this case cannot occur. If $c_1^2(M) = 2\chi(M) +
3\sigma(M) \geq 0$, this
leads to a contradiction, taking into account relation (\ref{temp00})
and the
fact that the signature has to be negative.
If $c_1^2(M) < 0$, first remark that $(M, \Om)$
must be a minimal symplectic 4-manifold, by the argument invoked above.
A. Liu \cite{Liu} shows that a minimal symplectic manifold with $b^+(M) =
1$ and
$c_1^2(M) < 0$, is an irrational ruled surface. This is a contradiction
since the signature of any ruled surface is 0, but
with our assumptions $\sigma(M) = 1 - b^-(M) < 0$.

\vspace{0.2cm}

Now we prove that $({\bar M}, g, {\bar J}, {\bar \Om})$ is, in fact, a
K{\"a}hler structure. Using the local considerations proved in Proposition
\ref{prop2} (ii), $({\bar M}, g, {\bar J}, {\bar \Om})$ is an almost \ka
structure and $(\na_X {\bar J})=0$ for any vector $X$ from ${\cal
D}^{\perp}$, where ${\cal D}$ is the \ka nullity of $(g, J)$. Let
${B, JB}$ be an orthonormal basis for the \ka nullity ${\cal D}$ of $J$
and
${A, JA}$ an orthonormal basis for ${\cal D}^{\perp}$.  Consider on $(M,
g, J)$
the {\it first canonical connection} $\na^0$ , defined by Lichnerowicz in
\cite{lichne} to be
$$\na^0_X Y=\na_X Y - \frac{1}{2}J (\na_X J)(Y).$$
This is a Hermitian connection, so its Ricci form $\gamma^0$ is, up to a
constant, a representative for the first Chern class of $(M, J)$. A short
computation of $\gamma^0$ gives for any almost Hermitian manifold (see
\cite{FH}):
$$ \gamma^0(X,Y) = 4Ric^*(JX,Y)  + < J\na_X J, \na_Y J> .$$
The condition $W^+_2=0$ is equivalent to the tensor $Ric^*$ being
symmetric (see (\ref{w^+_2})). 
Since we also assume that the Ricci tensor is $J$-invariant, it then follows
that $Ric^*$ and $Ric$ have common traceless part (see (\ref{U2})).
Using now Proposition \ref{prop2},(i), the above expression for $\gamma^0$ 
simplifies in
our case to:
$$ \gamma^0 = (s + \kappa) A \wedge JA + (s^* - \kappa) B \wedge JB .$$
Therefore,
\begin{equation}\label{temp2}
c_1^2(M)=\frac{1}{32\pi^2}\int_M (s + \kappa)(s^* - \kappa) \; dV_g.
\end{equation}
Concerning $({\bar M},g,{\bar J})$, we obtain similarly
$$ {\bar \gamma}^0 = - ({\bar s}^* + \kappa) A \wedge JA + (s - \kappa) B
\wedge JB ,$$
hence
\begin{equation}\label{temp3}
{\bar c}_1^2({\bar M})= \frac{1}{32\pi^2}\int_M ({\bar s}^* +
\kappa)(\kappa
- s) \;
dV_g,
\end{equation}
where ${\bar s}^*$ is the star-scalar curvature of $(g,{\bar J})$.
Taking into account that $c_1^2(M)={\bar c}_1^2({\bar M})=0$ and
$$\kappa - s =  - 3(s^* - \kappa) = \frac{3}{4} |\na J|^2 = const >0, $$
we
obtain
from (\ref{temp2}) and (\ref{temp3})
$\int_M ({\bar {\kappa}} - s) \; dV_g=0.$
Thus, $(g,{\bar J})$ is K{\"a}hler according to (\ref{s-s^*}).\\
Moreover, from (\ref{temp2}) we get 
$\int_M (s +\kappa) \; dV_g = 0$, 
which implies $\int_M s \; dV_g < 0$ (see (\ref{s-s^*})). It is well
known
that the anti-canonical bundle of a compact K{\"a}hler surface of negative
total scalar curvature has no holomorphic sections; then $(M,{\bar
J})$ is either a ruled surface with base of genus at 
least 2, or else
the Kodaira dimension is at least 1 (see for example \cite{BPV}). If
$({\bar M},{\bar J})$ is a ruled 
surface, then it is minimal, since $\sigma({\bar M})=0$, and hence
$\chi({\bar M})<0$, a contradiction. Moreover, the Kodaira dimension
of  $({\bar M},{\bar J})$ could not be 2 since any  surface of general
type has positive Euler 
number. We conclude that $({\bar M},{\bar J})$ belongs to class VI.  $\square$

\section{Pointwise constant totally-real sectional curvature}

Let $(M, g, J, \Om)$ be an almost Hermitian manifold.  As we have
already mentioned in the introduction, a two-plane 
$\sigma= X \wedge Y$ in the tangent bundle $TM$ is said to be 
{\it totally-real}
if $J \sigma = JX \wedge JY$ is orthogonal to $\sigma$. The
almost Hermitian manifold has {\it pointwise constant totally-real
sectional curvature} if at any point the sectional
curvature of the metric $g$ is the same on all totally-real planes at that
point.

The condition that $\si = X \wedge Y$ is a totally-real
plane is clearly equivalent to $g(JX, Y) = 0$, i.e., to $\si$ being a
Lagrangian plane with respect to the fundamental 2-form $\Om$. This
allows us to obtain the following simple characterization of the
4-dimensional
almost Hermitian manifolds of pointwise constant
totally-real sectional curvature in terms of the $U(2)$-decomposition
of the curvature tensor:
\begin{Lemma}\label{lem6} {\rm(\cite{Gan},\cite{AGI1})}
An almost Hermitian 4-manifold $(M, g, J)$ is of pointwise constant
totally-real sectional curvature
if and only if the Ricci tensor is $J$-invariant, $W_3 ^+ = 0$, $W^- = 0$.
The totally-real sectional curvature $\miu$ is given in this case by
$\miu = \frac{2s - \kappa}{24}$.
\end{Lemma}
{\it Proof:} At any point $x \in M$ we denote by ${\cal J}_x ^{\perp}$ the
space of all positive Hermitian structures of $(T_x M, g)$, anti-commuting
with $J$. The space  ${\cal J}_x ^{\perp}$  can be identified with the
elements of $\La_x ^{anti} M$ of square norm 2, via the metric $g$. We
start with the observation that a 2-plane $\si = X\wedge Y$ in $T_x M$ is
totally-real with respect to $J$ if and only if there exists
$ I \in {\cal J}_x ^{\perp}$ such that $\si $ is a {\it holomorphic} plane with
respect to 
$I$. Hence,  the set of holomorphic planes
for $I$  when $I$ varies in ${\cal J}_x ^{\perp}$
coincide with the set of totally-real planes for $J$.
Thus, $(M, g, J)$ is of pointwise constant totally-real sectional
curvature if and only if for any
$ I \in {\cal J}_x ^{\perp}$, $(M, g, I)$ is of pointwise constant
holomorphic sectional curvature,  say $c$, which does not depend
on $I$. In particular we have that $W^- = 0$ and $ Ric_0 $ is
$I$-anti-invariant for any $ I \in {\cal J}_x ^{\perp}$, cf.
\cite{kod,ADM}.
The second condition
is equivalent to the Ricci tensor (equivalently, $Ric_0$) being
$J$-invariant. Finally computing the holomorphic sectional curvatures
for $I$ and $K = J\circ I$, we obtain that they are equal if and only if
$W_{3}^+ = 0$, i.e., if and only if  $R$ has the form:
$$ R = \frac{s}{12} id + {\tilde {Ric_0}}^{inv} + \frac{\kappa}{8}
\Om \otimes \Om - \frac{\kappa}{12} id^{+} .$$
Then for any totally-real 2-plane $\sigma$ we have
$$ <R (\sigma), \sigma> = \frac{s}{12} |\sigma|^2 - \frac{\kappa}{12}
|\sigma_{+}|^2 = \frac{(2s - \kappa)}{24} |\sigma|^2 ,$$
which shows that the totally-real sectional curvature
is given by the formula claimed. $\square$

\vspace{0.1cm}
\noindent 
{\bf Proof of Theorem 2:} 
By Lemma \ref{lem6} the Ricci tensor is $J$-invariant, so
in the case on non-negative scalar curvature we apply
\cite[Theorem 1]{Dr1}. Suppose now that the scalar curvature is
negative.
 Since $g$ is self-dual (see Lemma \ref{lem6}), the Bach tensor vanishes (see
(\ref{gauduchon}) and
it follows from Corollary \ref{c2} that $g$ is Einstein. Thus, $(M,g,J,\Om)$ 
is a But compact almost
\ka Einstein mainifolds for which $W^+_3=0$ (see Lemma \ref{lem6}). It must be
K{\"a}hler according to \cite[Theorem 4.3.4]{Ar}. $\square$

\vspace{0.2cm}

\noindent
{\bf Remark 5.} Theorem 2 can be slightly generalized by
assuming that the scalar curvature is everywhere non-negative
or non-positive (but not necessarily constant). Indeed, in the non-negative
case,  the integrability of the
almost \ka structure follows by \cite{Dr1}. In the non-positive
case by Remark 2 we have that either the scalar
curvature identically vanishes, or else the metric is Einstein. If the
scalar curvature identically vanishes, we apply again \cite{Dr1}, while
in the latter case the the integrability of the almost \ka structure
follows from \cite{Ar}. 
We also observe that every compact almost \ka 4-manifold of pointwise
constant non-negative totally-real sectional curvature is K{\"a}hler. Indeed,
let $\mu$ be  the totally-real sectional curvature.
By (\ref{s-s^*}) and (\ref{kappa}), we have at any point,
$ s \leq \kappa$. In particular, we obtain from Lemma \ref{lem6} the
following pointwise
inequality for the totally-real sectional curvature: $\miu \leq
\frac{s}{24}$.
Hence, $\miu \geq 0$ implies $s\geq 0$ at any point and the result
follows from the observations above. $\Diamond$

\vspace{0.2cm}

 As pointed out by Derdzinski \cite{De}, any \ka metric in
dimension 4 has a degenerate spectrum for the positive Weyl tensor, that
is, the endomorphism $W^+$
of $\Lambda^+M$ has, at each point, at most two distinct
eigen-values. Derdzinski also remarked that all known examples of
Einstein metrics on compact, orientable 4-manifolds have degenerate
spectrum for $W^+$ (with one of the orientations of the manifold, at
least). To our knowledge, this remark is still valid at this date.  
It is not hard to see that any compatible almost \ka structure with 
an Einstein metric with degenerate spectrum for $W^+$, on a 
compact, {\it oriented}, 4-manifold, has to be, in fact,
K{\"a}hler. (In other words, the Goldberg conjecture is true for such
metrics.) Indeed, this follows essentially from the classification result
of
Derdzinski \cite[Theorem 2]{De}, with touches of some more recent works in
each of the three cases that occur:
\newline
(i) \cite[Theorem 2.4]{Ar3} or \cite[Theorem 1]{OS2} for the case of
Einstein anti-self-dual metric;
\newline
(ii) \cite[Theorem 1]{ApDr} for the case when the metric admits a \ka
structure;
\newline
(iii) \cite{Se2} in the last case, when the scalar curvature is positive.

\vspace{0.2cm}

Regarding the degeneracy of the
spectrum
of the positive Weyl tensor for an almost \ka
4-manifold
of pointwise constant totally-real sectional curvature, we have the
following
\begin{theo}\label{th4}
Let $(M, g, J, \Om)$ be a compact almost \ka 4-manifold
of pointwise constant totally-real sectional curvature. Then $J$ is
integrable if and only if the spectrum of the positive
Weyl tensor is degenerate.
\end{theo}
{\it Proof:} It is known that for an almost Hermitian 4-manifold with $W_3^+
= 0$,
the spectrum of $W^+$ is degenerate if and only if $W_2^+ = 0$ (see
\cite{Sa,AG}).
Thus, according to Lemmas \ref{lem5} and \ref{lem6}, for an almost Hermitian
4-manifold of pointwise constant totally-real sectional curvature, the
spectrum of
$W^+$ is degenerate if and only if $(g, J)$ satisfies the second
curvature condition of Gray. If $(M,g,J)$ 
is almost \ka {\it non-K{\"a}hler} 4-manifold of pointwise constant
totally-real sectional curvature, whose positive Weyl
tensor is
degenerate, it follows from Proposition 2 that there is an almost K{\"a}hler
structure ${\bar J}$, compatible with $g$ and the inverse orientation of
$M$, such that the Ricci tensor of $g$ is ${\bar J}$-invariant. As $g$ is
self-dual (Lemma \ref{lem6}), it follows from Corollary \ref{c1} that
$(g,{\bar J})$ is K{\"a}hler, hence $g$ is scalar flat, contradiction to
\cite[Theorem 1]{Dr1} $\square$

\end{document}